# A Generalization of Certain Remarkable Points of the Triangle Geometry

Prof. Claudiu Coandă – National College "Carol I", Craiova, Romania
Prof. Florentin Smarandache – University of New Mexico, Gallup, U.S.A.
Prof. Ion Pătrașcu – National College "Frații Buzești", Craiova, Romania

In this article we prove a theorem that will generalize the concurrence theorems that are leading to the Franke's point, Kariya's point, and to other remarkable points from the triangle geometry.

**Theorem 1:**

Let $P(\alpha, \beta, \gamma)$ and $A', B', C'$ its projections on the sides $BC$, $CA$ respectively $AB$ of the triangle $ABC$.

We consider the points $A''$, $B''$, $C''$ such that $\overrightarrow{PA''} = k\overrightarrow{PA'}$, $\overrightarrow{PB''} = k\overrightarrow{PB'}$, $\overrightarrow{PC''} = k\overrightarrow{PC'}$, where $k \in R^*$. Also we suppose that $AA', BB', CC'$ are concurrent. Then the lines $AA'', BB'', CC''$ are concurrent if and only if are satisfied simultaneously the following conditions:

$$\alpha\beta c\left(\frac{\beta}{b}\cos A - \frac{\alpha}{a}\cos B\right) + \beta\gamma a\left(\frac{\gamma}{c}\cos B - \frac{\beta}{b}\cos C\right) + \gamma\alpha b\left(\frac{\alpha}{a}\cos C - \frac{\gamma}{c}\cos A\right) = 0$$

$$\frac{\alpha^2}{a^2}\cos A\left(\frac{\gamma}{c}\cos B - \frac{\beta}{b}\cos C\right) + \frac{\beta^2}{b^2}\cos B\left(\frac{\alpha}{a}\cos C - \frac{\gamma}{c}\cos A\right) + \frac{\gamma^2}{c^2}\cos C\left(\frac{\beta}{b}\cos A - \frac{\alpha}{a}\cos B\right) = 0$$

**Proof:**

We find that

$$A'\left(0, \frac{\alpha}{2a^2}(a^2+b^2-c^2)+\beta, \frac{\alpha}{2a^2}(a^2-b^2+c^2)+\gamma\right)$$

$$\overrightarrow{PA''} = k\overrightarrow{PA'} = k\left[-\alpha\overrightarrow{r_A} + \frac{\alpha}{2a^2}(a^2+b^2-c^2)\overrightarrow{r_B} + \frac{\alpha}{2a^2}(a^2-b^2+c^2)\overrightarrow{r_C}\right]$$

$$\overrightarrow{PA''} = (\alpha''-\alpha)\overrightarrow{r_A} + (\beta''-\beta)\overrightarrow{r_B} + (\gamma''-\gamma)\overrightarrow{r_C}$$

We have:

$$\begin{cases} \alpha''-\alpha = -k\alpha \\ \beta''-\beta = \dfrac{k\alpha}{2a^2}(a^2+b^2-c^2), \\ \gamma''-\gamma = \dfrac{k\alpha}{2a^2}(a^2-b^2+c^2) \end{cases}$$

Therefore:

$$\begin{cases} \alpha'' = (1-k)\alpha \\ \beta'' = \dfrac{k\alpha}{2a^2}(a^2+b^2-c^2) + \beta \\ \gamma'' = \dfrac{k\alpha}{2a^2}(a^2-b^2+c^2) + \gamma \end{cases}$$

Hence:
$$A''\left((1-k)\alpha,\ \dfrac{k\alpha}{2a^2}(a^2+b^2-c^2)+\beta,\ \dfrac{k\alpha}{2a^2}(a^2-b^2+c^2)+\gamma\right)$$

Similarly:
$$B'\left(-\dfrac{\beta}{2b^2}(-a^2-b^2+c^2)+\alpha,\ 0,\ -\dfrac{\beta}{2b^2}(a^2-b^2-c^2)+\gamma\right)$$

$$B''\left(-\dfrac{k\beta}{2b^2}(-a^2-b^2+c^2)+\alpha,\ (1-k)\beta,\ -\dfrac{k\beta}{2b^2}(a^2-b^2-c^2)+\gamma\right)$$

$$C'\left(-\dfrac{\gamma}{2c^2}(-a^2+b^2-c^2)+\alpha,\ -\dfrac{\gamma}{2c^2}(a^2-b^2-c^2)+\beta,\ 0\right)$$

$$C''\left(-\dfrac{k\gamma}{2c^2}(-a^2+b^2-c^2)+\alpha,\ -\dfrac{k\gamma}{2c^2}(a^2-b^2-c^2)+\beta,\ (1-k)\gamma\right)$$

Because $AA'$, $BB'$, $CC'$ are concurrent, we have:

$$\dfrac{-\dfrac{\alpha}{2a^2}(-a^2-b^2+c^2)+\beta}{-\dfrac{\alpha}{2a^2}(-a^2+b^2-c^2)+\gamma} \cdot \dfrac{-\dfrac{\beta}{2b^2}(-a^2-b^2-c^2)+\gamma}{-\dfrac{\beta}{2b^2}(-a^2+b^2+c^2)+\alpha} \cdot \dfrac{-\dfrac{\gamma}{2c^2}(-a^2+b^2-c^2)+\alpha}{-\dfrac{\gamma}{2c^2}(a^2-b^2-c^2)+\beta} = 1$$

We note
$$M = \dfrac{\alpha}{2a^2}(a^2+b^2-c^2) = \dfrac{\alpha}{a}\cdot b\cos C$$

$$N = \dfrac{\alpha}{2a^2}(a^2-b^2+c^2) = \dfrac{\alpha}{a}\cdot c\cos B$$

$$P = \dfrac{\beta}{2b^2}(-a^2+b^2+c^2) = \dfrac{\beta}{b}\cdot c\cos A$$

$$Q = \dfrac{\beta}{2b^2}(a^2+b^2-c^2) = \dfrac{\beta}{b}\cdot a\cos C$$

$$R = \dfrac{\gamma}{2c^2}(a^2-b^2+c^2) = \dfrac{\gamma}{c}\cdot a\cos B$$

$$S = \dfrac{\gamma}{2c^2}(-a^2+b^2+c^2) = \dfrac{\gamma}{c}\cdot a\cos A$$

The precedent relation becomes
$$\dfrac{M+\beta}{N+\gamma}\cdot\dfrac{P+\gamma}{Q+\alpha}\cdot\dfrac{R+\alpha}{S+\beta} = 1$$

The coefficients $M, N, P, Q, R, S$ verify the following relations:

$$M + N = \alpha$$
$$P + Q = \beta$$
$$R + S = \gamma$$

$$\frac{M}{Q} = \frac{\alpha}{\beta} \cdot \frac{b^2}{a^2} = \frac{\frac{\alpha}{a^2}}{\frac{\beta}{b^2}}$$

$$\frac{P}{S} = \frac{\beta}{\gamma} \cdot \frac{c^2}{b^2} = \frac{\frac{\beta}{b^2}}{\frac{\gamma}{c^2}}$$

$$\frac{R}{N} = \frac{\gamma}{\alpha} \cdot \frac{a^2}{c^2} = \frac{\frac{\gamma}{c^2}}{\frac{\alpha}{a^2}}$$

Therefore $\frac{M}{Q} \cdot \frac{P}{S} \cdot \frac{R}{N} = 1$

$$(M+\beta)(P+\gamma)(R+\alpha) = \alpha\beta\gamma + \alpha\beta P + \beta\gamma R + \gamma\alpha M + \alpha MP + \beta PR + \gamma RM + MPR$$
$$(N+\gamma)(Q+\alpha)(S+\beta) = \alpha\beta\gamma + \alpha\beta N + \beta\gamma Q + \gamma\alpha S + \alpha NS + \beta NQ + \gamma QS + NQS.$$

We deduct that:
$$\alpha\beta P + \beta\gamma R + \gamma\alpha M + \alpha MP + \beta PR + \gamma RM = \alpha\beta N + \beta\gamma Q + \gamma\alpha S + \alpha NS + \beta NQ + \gamma QS + NQS \quad (1)$$

We apply the theorem:

Given the points $Q_i(a_i, b_i, c_i)$, $i = \overline{1,3}$ in the plane of the triangle $ABC$, the lines $AQ_1, BQ_2, CQ_3$ are concurrent if and only if $\frac{b_1}{c_1} \cdot \frac{c_2}{a_2} \cdot \frac{a_3}{b_3} = 1$.

For the lines $AA", BB", CC"$ we obtain
$$\frac{kM+\beta}{kN+\gamma} \cdot \frac{kP+\alpha}{kS+\beta} \cdot \frac{kR+\alpha}{kS+\beta} = 1.$$

It result that
$$k^2(\alpha\beta P + \beta\gamma R + \gamma\alpha M) + k(\alpha MP + \beta PR + \gamma RM) =$$
$$= k^2(\alpha\beta N + \beta\gamma Q + \gamma\alpha S) + k(\alpha NS + \beta NQ + \gamma QS) \quad (2)$$

For relation (1) to imply relation (2) it is necessary that
$$\alpha\beta P + \beta\gamma R + \gamma\alpha M = \alpha\beta N + \beta\gamma Q + \gamma\alpha S$$

and
$$\alpha NS + \beta NQ + \gamma QS = \alpha MP + \beta PR + \gamma RM$$

or

$$\begin{cases} \alpha\beta c\left(\dfrac{\beta}{b}\cos A - \dfrac{\alpha}{a}\cos B\right) + \beta\gamma a\left(\dfrac{\gamma}{c}\cos B - \dfrac{\beta}{b}\cos C\right) + \gamma\alpha b\left(\dfrac{\alpha}{a}\cos C - \dfrac{\gamma}{c}\cos A\right) = 0 \\ \dfrac{\alpha^2}{a^2}\cos A\left(\dfrac{\gamma}{c}\cos B - \dfrac{\beta}{b}\cos C\right) + \dfrac{\beta^2}{b^2}\cos B\left(\dfrac{\gamma}{c}\cos B - \dfrac{\beta}{b}\cos C\right) + \dfrac{\gamma^2}{c^2}\cos C\left(\dfrac{\beta}{b}\cos A - \dfrac{\alpha}{a}\cos B\right) = 0 \end{cases}$$

As an open problem, we need to determine the set of the points from the plane of the triangle $ABC$ that verify the precedent relations.

We will show that the points $I$ and $O$ verify these relations, proving two theorems that lead to Kariya's point and Franke's point.

**Theorem 2** (Kariya -1904)

Let $I$ be the center of the circumscribe circle to triangle $ABC$ and $A', B', C'$ its projections on the sides $BC$, $CA$, $AB$. We consider the points $A'', B'', C''$ such that:
$$\overrightarrow{IA''} = k\overrightarrow{IA'}, \ \overrightarrow{IB''} = k\overrightarrow{IB'}, \ \overrightarrow{IC''} = k\overrightarrow{IC'}, \ k \in R^*.$$
Then $AA''$, $BB''$, $CC''$ are concurrent (the Kariya's point)

**Proof:**

The barycentric coordinates of the point $I$ are $I\left(\dfrac{a}{2p}, \dfrac{b}{2p}, \dfrac{c}{2p}\right)$.

Evidently:
$$abc(\cos A - \cos B) + abc(\cos B - \cos C) + abc(\cos C - \cos A) = 0$$
and
$$\cos A(\cos B - \cos C) + \cos B(\cos C - \cos A) + \cos C(\cos A - \cos B) = 0.$$

In conclusion $AA''$, $BB''$, $CC''$ are concurrent.

**Theorem 3** (de Boutin - 1890)

Let $O$ be the center of the circumscribed circle to the triangle $ABC$ and $A', B', C'$ its projections on the sides $BC$, $CA$, $AB$. Consider the points $A'', B'', C''$ such that $\dfrac{OA'}{OA''} = \dfrac{OB'}{OB''} = \dfrac{OC'}{OC''} = k, \ k \in R^*$. Then the lines $AA''$, $BB''$, $CC''$ are concurrent (The point of Franke – 1904).

**Proof:**

$O\left(\dfrac{R^2}{2S}\sin 2A, \dfrac{R^2}{2S}\sin 2B, \dfrac{R^2}{2S}\sin 2C\right)$, $P = N$, because $\dfrac{\sin 2B \cos A}{\sin B} - \dfrac{\sin 2A \cos B}{\sin A} = 0$.

Similarly we find that $R = Q$ and $M = S$.

Also $\alpha MP = \alpha NS$, $\beta PR = \beta NQ$, $\gamma RM = \gamma QS$. It is also verified the second relation from the theorem hypothesis. Therefore the lines $AA''$, $BB''$, $CC''$ are concurrent in a point called the Franke's point.

**Remark 1**:

It is possible to prove that the Franke's points belong to Euler's line of the triangle $ABC$.

**Theorem 4**:

Let $I_a$ be the center of the circumscribed circle to the triangle $ABC$ (tangent to the side $BC$) and $A'$, $B'$, $C'$ its projections on the sites $BC$, $CA$, $AB$. We consider the points $A"$, $B"$, $C"$ such that $\overrightarrow{IA"} = k\overrightarrow{IA'}$, $\overrightarrow{IB"} = k\overrightarrow{IB'}$, $\overrightarrow{IC"} = k\overrightarrow{IC'}$, $k \in R^*$. Then the lines $AA"$, $BB"$, $CC"$ are concurrent.

**Proof**

$$I_a\left(\frac{-a}{2(p-a)}, \frac{b}{2(p-a)}, \frac{c}{2(p-a)}\right);$$

The first condition becomes:
$$-abc(\cos A + \cos B) + abc(\cos B - \cos C) - abc(-\cos C - \cos A) = 0,$$
and the second condition:
$$\cos A(\cos B - \cos C) + \cos B(-\cos C - \cos A) + \cos C(\cos A + \cos B) = 0$$
Is also verified.

From this theorem it results that the lines $AA"$, $BB"$, $CC"$ are concurrent.

**Observation 1:**

Similarly, this theorem is proven for the case of $I_b$ and $I_c$ as centers of the ex-inscribed circles.